\documentclass[reqno,12pt]{amsart}
\usepackage{amsmath, latexsym, amsfonts, amssymb, amsthm, amscd}

\numberwithin{equation}{section}

\newtheorem{theorem}{Theorem}[section]
\newtheorem{lemma}[theorem]{Lemma}
\newtheorem{proposition}[theorem]{Proposition}
\newtheorem{corollary}[theorem]{Corollary}
\newtheorem{rem}[theorem]{Remark}

\renewcommand{\tilde}{\widetilde}          
\DeclareMathSymbol{\leqslant}{\mathalpha}{AMSa}{"36} 
\DeclareMathSymbol{\geqslant}{\mathalpha}{AMSa}{"3E} 
\DeclareMathSymbol{\eset}{\mathalpha}{AMSb}{"3F}     
\renewcommand{\leq}{\;\leqslant\;}                   
\renewcommand{\geq}{\;\geqslant\;}                   

\newcommand{\R}{\mathbb{R}}
\newcommand{\Z}{\mathbb{Z}}
\newcommand{\N}{\mathbb{N}}

\def\Pd{{\mathtt P}}
\def\Ed{{\mathtt E}}
\def\8{{\infty}}

\title[Majorizing cascades for directed polymers 
]{Majorizing multiplicative cascades for directed polymers 
in random media}\thanks{Draft \today. 
}
\author{Francis Comets, Vincent Vargas}
\thanks{Partially
   supported by CNRS (UMR 7599
``Probabilit{\'e}s et Mod{\`e}les
Al{\'e}atoires'')}

\begin{document}

\maketitle
\begin{center}
{\footnotesize \noindent
 Universit{\'e} Paris 7,\\
Math{\'e}matiques, case 7012,\\ 2, place Jussieu, 75251 Paris, France}

{\footnotesize \noindent e-mail: \texttt{comets@math.jussieu.fr,vargas@math.jussieu.fr}}
\end{center}

\begin{abstract}
In this note we give upper bounds for the free energy of discrete
directed polymers in random media. The bounds are given by the so-called
generalized multiplicative 
cascades from the statistical theory of turbulence. For the polymer
model, we derive that the quenched free energy is different from
the annealed one in dimension 1, for any finite temperature and
general  environment. This implies localization of the polymer.
\end{abstract}
\vspace{1cm}
\footnotesize
\noindent{\bf Short Title.} Majorizing cascades for directed polymers

\noindent{\bf Key words and phrases.} Directed polymers, random
environment, strong disorder, generalized multiplicative 
cascades.

\noindent{\bf MSC 2000 subject classifications.} Primary 60K35;
secondary 60J30, 82D30, 82A51

\normalsize

\section{Introduction}
Let  $\omega=(\omega_{n})_{n \in \N}$ be  the simple
random walk on the $d$-dimensional integer lattice $\Z^d$ starting at
0, defined on
a probability space $(\Omega, \mathcal{F}, P)$.
We also consider a sequence $\eta=(\eta(n,x))_{(n,x) \in \N
  \times \Z^d}$ of real valued, non-constant and
i.i.d. random variables defined on another probability space $(H,
\mathcal{G}, Q)$ with finite exponential moments. 
The path $\omega$ represents the directed polymer and 
$\eta$ the random environment.

For any $n>0$, we define the (random) polymer measure $\mu_{n}$ on
the path space $(\Omega, \mathcal{F})$ by:
$$\mu_{n}(d\omega)=\frac{1}{Z_{n}}\exp(\beta H_{n}(\omega))P(d\omega)$$
where $\beta \in \R^+$ is the inverse temperature, where
$$
 H_{n}(\omega)\overset{\rm def.}{=}\sum_{j=1}^{n}\eta(j,\omega_{j}) \qquad
$$
and where
$$
Z_{n}=P[\exp(\beta H_{n}(\omega))] \qquad
$$
is the partition function. We use the notation $P[X]$ for the
expectation of a random variable $X$.
By symmetry, we can -- and we will -- restrict to $\beta \geq 0$. 
\medskip

 The free energy of the polymer is  defined as the limit 
\begin{equation} \label{def:p}
p(\beta) =\lim_{n \to \infty}\frac{1}{n}\ln(Z_{n}(\beta)/Q[Z_{n}(\beta)]) 
\end{equation}
where the limit exists $Q$-a.s. and in $L^p$ for all $p
\geq 1$ and is constant (cf. \cite{cf: CoShYo}).
An application of Jensen's inequality to the concave function $\ln (\cdot)$
yields $p(\beta) \leq 0$. 
As shown in theorem 3.2  (b) in
\cite{cf: CoYo},
there exists a $\beta_{c} \in [0, \infty]$ such that 
\begin{equation*}
p(\beta) \quad
\begin{cases}
 \quad = 0& \qquad \text{if} \quad \beta \in [0,\beta_{c}], \\
 \quad < 0& \qquad \text{if} \quad \beta > \beta_{c}.
\end{cases}
\end{equation*}
An
important question in the study of directed polymers is to find the
$\beta$ such that $p(\beta) < 0$. Indeed, one can show that the
negativity of $p(\beta)$ is equivalent to a localization property for
$(\omega_{n})_{n \in \N}$,$(\tilde{\omega}_{n})_{n \in \N}$ two
independent random walks under the polymer measure $\mu_{n}$ (cf
Corollary 2.2 in \cite{cf: CoShYo}):
\begin{equation*}
p(\beta) < 0 \qquad \Longleftrightarrow \qquad \exists c > 0 \quad
\varliminf_{n \to \infty}\frac{1}{n}\sum_{k=1}^{n}\mu_{k-1}^{\otimes 2}(\omega_{k}=\tilde{\omega}_{k})\geq c \qquad Q.a.s.
\end{equation*}
The statement in the right-hand side means that the polymer localizes
in narrow corridors with positive probability. 
It is not known how to characterize  directly these corridors, and therefore
this criterion for the transition localization/delocalization 
is rather efficient since it does not require any knowledge on them.
 Hence, it is important to get good upper bounds 
on $p$ in order to spot the transition. Our main result is the
following.

\begin{theorem} \label{mainth}
In dimension $d=1$, $\beta_{c}=0$.
\end{theorem}

There is a clear consensus on this fact in the physics
literature, but no proof for it, except via the replica method
or in the (different) case of a space-periodic environment 
where much more computations can be
performed
\cite{BrDe00}.

This result follows from a family of upper bounds, given by the free
energies of models on trees depending on an integer parameter $m$
($m \geq 1$).
These trees are deterministic and regular, with random weights, they
fall in the scope of the generalized multiplicative
cascades \cite{cf:Liu} or smoothing transformations \cite{cf:DuLi83}
which are well known generalizations of the random cascades introduced 
in \cite{cf:Man} for a statistical description of turbulence.
When the environment variables have nice concentration properties
-- e.g., gaussian or bounded $\eta$'s --, 
we prove  in theorem \ref{th:varfor}
that the polymer free energy is the infimum over $m$ of the one of
the $m$-tree  model. For general environmental distribution we only
have an upper bound from theorem \ref{lem: ine}, but it is enough to show
the above theorem.
This also explains the title of the present paper.

Recall at this point that  directed polymers in a Bernoulli 
random environment are positive temperature versions of oriented percolation.
Our bounds here have a flavor similar to the lower bounds for the
critical threshold in 2-dimensional oriented percolation (i.e., $d=1$
in our notations)
in section 6 of Durrett \cite{cf:Dur84}. In that paper, percolation is
compared to Galton-Watson processes obtained in running oriented percolation
for $m$ steps ($m \geq 1$), and then using the distribution of wet
sites as offspring distribution. 

Next, we comment on the case of supercritical 1-dimensional oriented 
percolation. Then, $\eta$ is Bernoulli distributed with parameter $p>\vec
p_c(1)$. The infinite cluster is the set of points $(t,x)$
with $t \in \N, x \in \Z, P(\omega_t=x)>0$, which are connected
to $\8$ by an open oriented path -- i.e., a path $\omega$ with
$\eta(s,\omega_s)=1\; \forall s \geq t$. It is known that this
cluster, at large scale, 
is approximatively a cone with vertex $(0,0)$, direction $[0,x)$ and
positive angle, and it has a positive density. In words, there is a
huge number of oriented paths of length $n$ with energy $H_n=n -
{\mathcal O}(1)$. However, according to the theorem,
 the polymer measure has a strong localization property. 
This first seems paradoxical, since there are exponentially many
suitable paths on the energetic level. Hence, this is essentially
an entropic phenomenon, due to large 
fluctuations in the number of such paths.

For numerics, our upper bounds
 do not seem very efficient: on the
basis of preliminary numerical simulations they converge quite slowly
as $m \to \8$.
Finally we mention that
lower bounds for  the polymer free energy can be obtained
from a well-known 
super-additivity property, see formula (\ref{eq: super+}).

\section{Notations and preliminaries}

We first introduce some further notations.

Let $((\omega_{n})_{n \in \N},(P^x)_{x \in \Z^d})$ denote the simple
random walk on the $d$-dimensional integer lattice $\Z^d$, defined on
a probability space $(\Omega, \mathcal{F})$: for $x$
in $\Z^d$, under the measure $P^x$, $(\omega_{n}-\omega_{n-1})_{n \geq
  1}$ are independent and 
\begin{equation*}
P^x(\omega_{0}=x)=1, \; \; P^x(\omega_{n}-\omega_{n-1}=\pm
\delta_{j})=\frac{1}{2d},\; \; j=1, \ldots, d, 
\end{equation*}
where $(\delta_{j})_{1 \leq j \leq d}$ is the j-th vector of the
canonical basis of $\Z^d$. Like in the introduction, we  will use the notation
$P$ for $P^0$. 

For the environment, we assume that  for all $\beta \in \R$,
$$
 \lambda(\beta)\overset{\rm def.}{=} \ln
Q(e^{\beta  \eta(n,x)})< \infty.
$$

It is convenient to consider the normalized partition function
$$
W_{n}= Z_{n}/Q[Z_n]
=P[\exp(\beta H_{n}(\omega)-n \lambda(\beta))].
$$
We define for $k < n$, $x,y \in \Z^d$,
\begin{equation*}
H_{k,n}(\omega)=\sum_{j=1}^{n-k}\eta(k+j,\omega_{j})
\end{equation*}
and 
\begin{equation} \label{eq:Wn(x)}
W_{k,n}^{x}(y)=P^{x}(e^{\beta H_{k,n}(\omega)-(n-k)\lambda(\beta)}1_{\omega_{n-k}=y}).
\end{equation}
In the sequel, $W_{n}(x)$ will stand for $W_{0,n}^{x}(0)$.
The Markov property of the simple random walk yields
\begin{equation}\label{eq: mark}
W_{n}=\sum_{x,y \in \Z^d}W_{k}(x)W_{k,n}^{x}(y).
\end{equation}
This identity will be extensively used in the sequel.

Finally, we recall (\cite{cf: CoShYo})
that with $p$ defined by (\ref{def:p}) it holds
\begin{equation}\label{eq: super+}
p(\beta)= \lim_{n \to \infty}\frac{1}{n}Q(\ln(W_{n}(\beta)))
 = \sup_{n \geq 1}\frac{1}{n}Q(\ln(W_{n}(\beta)))  
\end{equation}
where the last equality is a consequence of super-additivity arguments.

\subsection{Definition and well known facts on generalized
  multiplicative cascades}
 In this section, we introduce a model of generalized multiplicative
 cascades on a tree. For an overview of results, we refer to \cite{cf:Liu}. Let $N \geq 2$ be a fixed integer and   
\begin{equation*}
U= \underset{k \in \N}{\bigcup}[\mid 1,N \mid]^{k}
\end{equation*}
be the set of all finite sequences $ u=u_{1} \ldots u_{k}$ of
elements in $[\mid 1,N \mid]$. With the previous notation, we write
$\mid u \mid =k $ for its length. For $u=u_{1} \ldots u_{k}$,$v=v_{1}
\ldots v_{k}$ two finite sequences, let $uv$ denote
the sequence $u_{1} \ldots u_{k}v_{1} \ldots v_{k}$. 
Let $q$ be a non degenerate probability
distribution on $(\R_{+}^{*})^{N}$. It is known (cf. \cite{cf:Liu})
that there exist a probability space 
with probability measure denoted by $\Pd$ (and expectation $\Ed$),
and random variables $(A_{u})_{u \in U}$ defined on this space, such that   
the random vectors $(A_{u1}, \ldots, A_{uN})_{u \in U}$ form an
i.i.d. sequence with common distribution $q$.      
We assume that 
the $(A_{i})_{1 \leq i \leq N}$ are normalized:
\begin{equation*}
\Ed(\sum_{i=1}^{N}A_{i})=1
\end{equation*}
and that they have moments of all order:
$ \Ed[\sum_{i=1}^{N}A_{i}^{p}] < \infty \;\forall p \in \R .$
Consider the  process
 $(W_{n}^{\rm casc})_{n \in \N}$  defined by  
\begin{equation} \label{def:casc}
W_{n}^{\rm casc}=\sum_{u_{1}, \ldots, u_{n} 
\in [\mid 1,N \mid]}A_{u_{1}}A_{u_{1}u_2}
\ldots
A_{u_{1}\ldots u_{n}}
\end{equation}
and the  filtration
\begin{equation*}
\mathcal{G}_{n}:= \sigma \lbrace A_{u}; \mid u \mid \leq n \rbrace,
\qquad n \geq 1.
\end{equation*}
Then $(W_{n}^{\rm casc},\mathcal{G}_{n})_{n \geq 1}$ 
is a non negative martingale
so the limit $W_{\infty}^{\rm casc}=\lim_{n \to \infty} W_{n}^{\rm
  casc}$ exists. We are
interested in the behavior of the associated free energy:
\begin{equation*} 
p_{n}=\frac{1}{n}\ln W_{n}^{\rm casc}\;.
\end{equation*}
In the case where the $(A_{i})_{i \leq N}$ are i.i.d, the exact
limit of $p_{n}$ as $n$ goes to infinty was derived in \cite{cf:
  Fra}. In the general case, the proofs in \cite{cf: Fra} can easily
be adapted to show the following summary result.
\begin{theorem} \label{th:p-casc}
The following convergence holds $\Pd$-a.s. and in $L^p$ for all $p \geq
1$:
\begin{equation*}
p_{n}\underset{n \to \infty}{\longrightarrow} 
\inf_{\theta \in
  ]0,1]}\frac{1}{\theta}\ln(\Ed\sum_{i=1}^{N}A_{i}^{\theta}) \leq 0,
\end{equation*}
\end{theorem}
where the inequality is a consequence of the normalization.
Finding the limit of $p_{n}$ as $n$ tends to infinity amounts to
studying the function $v$ defined by 
\begin{equation*}  
\forall \theta \in ]0,1], \qquad v(\theta)=\frac{1}{\theta}
\ln(\Ed\sum_{i=1}^{N}A_{i}^{\theta})
\;,
\end{equation*}
which has derivative 
$$
 v'(1)= \Ed\sum_{i=1}^{N}A_{i}\ln(A_{i})\;.$$
\begin{lemma}\label{lem: v}
If $\Ed\sum_{i=1}^{N}A_{i}\ln(A_{i}) \leq 0$, the function $v$ is
strictly decreasing on $]0,1]$ and thus 
\begin{equation*}
\inf_{\theta \in ]0,1]}v(\theta)=v(1)=0.
\end{equation*}
If $\Ed\sum_{i=1}^{N}A_{i}\ln(A_{i})>0$, there exists a unique
$\theta^{*} \in ]0,1[$ such that 
\begin{equation*}
\inf_{\theta \in ]0,1]}v(\theta)=v(\theta^{*})<0.
\end{equation*}
\end{lemma}
\proof
For all $\theta \in ]0,1]$, we have the following expression for the
derivative of $v$:
\begin{equation*}
v'(\theta)=\frac{g(\theta)}{\theta^2}
\end{equation*}
where $g$ is given by 
\begin{equation*}
g(\theta)=\theta\frac{\Ed\sum_{i=1}^{N}A_{i}^{\theta}\ln(A_{i})}{\Ed\sum_{i=1}^{N}A_{i}^{\theta}}-\ln(\Ed\sum_{i=1}^{N}A_{i}^{\theta}).
\end{equation*}
In particular, we obtain the value of $v'(1)$ given above.
By direct computation, one can obtain the following expression for $g'$  
\begin{equation*}
\forall \theta > 0 \qquad g'(\theta)= \theta\frac{\Ed(\sum_{i=1}^{N}A_{i}^{\theta}(\ln(A_{i})-\Ed(\ln(A)\mid A^{\theta}))^2)}{\Ed(\sum_{i=1}^{N}A_{i}^{\theta})}
\end{equation*}
where $\Ed(\ln(A)\mid A^{\theta})$ is a notation for 
\begin{equation*}
\Ed(\ln(A)\mid A^{\theta})=\frac{\Ed(\sum_{i=1}^{N}A_{i}^{\theta}\ln(A_{i}))}{\Ed(\sum_{i=1}^{N}A_{i}^{\theta})}.
\end{equation*}
In particular, $g$ is strictly increasing and we have 
\begin{equation*}
g(1)=\Ed(\sum_{i=1}^{N}A_{i}\ln(A_{i})).
\end{equation*}
By considering the two cases $g(1) \leq 0$ and $g(1) > 0$, we can
easily conclude. 
\qed

\subsection{Concentration of measure in the gaussian and the bounded
  case} 

For a complete survey on the concentration of measure phenomenon, we
refer to \cite{cf: Led}. In the gaussian case, we have
\begin{theorem}
Let $M \geq 1$ be an integer. We consider $\R^{M}$ equiped with the
usual euclidian norm $\|\cdot \|$. If $X_{M}$ is a standard gaussian
vector on some probability space (with a  probability measure $\Pd$)
 and $F$ is a
$C$-lipschitzian function ($|F(x)-F(y)| \leq C \|x-y\|$)
from $\R^{M}$ to $\R$ then
\begin{equation}\label{eq: con}    
\Ed(e^{\lambda(F(X_{M})-\Ed(F(X_{M})))}) \leq e^{\frac{C^2 \lambda^{2}}{2}}.
\end{equation}
Therefore, we have the following concentration result
\begin{equation}\label{eq: con1}
\Pd(\mid F(X_{M})-\Ed(F(X_{M})) \mid \geq r) \leq 2 e^{-\frac{r^2}{2C^{2}}} 
\end{equation}
\end{theorem}

In the bounded case, we get a similar concentration result
(cf. Corollary 3.3 in \cite{cf: Led}) . 
\begin{theorem}
Let $M \geq 1$ be an integer and $a<b$ be two real numbers. 
If $X_{M}$ is a random 
vector in $[a,b]^{M}$ with i.i.d. components  on some probability
space 
and $F$ is a convex and
$C$-lipschitzian function from $[a,b]^{M}$ to $\R$ for the euclidian norm,
then
\begin{equation}\label{eq: con'}    
\Ed(e^{\lambda(F(X_{M})-\Ed(F(X_{M})))}) \leq e^{C^2 (b-a)^2 \lambda^{2}}.
\end{equation}
Therefore, we have the following concentration result
\begin{equation}\label{eq: con1'}
\Pd( F(X_{M})-\Ed(F(X_{M})) \geq r) \leq  e^{-\frac{r^2}{4C^{2}(b-a)^2}} 
\end{equation}
\end{theorem}

We can derive from the above theorems a concentration result for the
free energy at time $n$: 
\begin{corollary}
If the environment $\eta$ is standard gaussian then for all $\lambda \geq
0$,
\begin{equation}\label{eq: conc}  
Q(e^{\lambda(\ln(W_{n})-Q(\ln(W_{n})))}) \leq 
e^{\frac{\beta^2  \lambda^2 n}{2}}. 
\end{equation}
If the environment $\eta$ belongs to $[a,b]$ for $a<b$ two real
numbers, then for all $\lambda \geq 0$,
\begin{equation}\label{eq: conc'}
Q(e^{\lambda(\ln(W_{n})-Q(\ln(W_{n})))}) \leq 
e^{\beta^2(b-a)^2 \lambda^2 n}.
\end{equation}
\end{corollary}
\proof
As a function of the environment, $\ln(W_{n})$ is convex and
$\beta\sqrt{n}$-lipschitzian (cf. the proof of proposition 1.4 in 
\cite{cf:CaHu02}). Therefore, in the gaussian case, the result is a
direct application of (\ref{eq: con}) and, in the bounded case, simply
  (\ref{eq: con'}).   
\qed

\section{Majorizing polymers with cascades}

Let us fix an integer $m \geq 1$ and define $L_{m}$ to be set of
points visited by the simple random walk at time $m$: 
\begin{equation*}
L_{m}\overset{def}{=}\lbrace x \in \Z^d; P(w_{m}=x)>0 \rbrace.
\end{equation*}
 We introduce $(W_{m,n}^{\rm tree})_{n \geq 1} \equiv
(W_{n}^{\rm casc})_{n \geq 1}$
 the martingale of the  
multiplicative cascade associated to the random vector 
$(W_{m}(x))_{x \in L_{m}}$, i.e., defined by (\ref{def:casc}) when 
$N=|L_m|$ and $q$
is the law of $(W_{m}(x))_{x \in L_{m}}$ with $W_{m}(x)$ from
(\ref{eq:Wn(x)}).
 Let $p_{m}^{ \rm tree}(\beta)$ denote the associated
free energy. In view of (\ref{th:p-casc}), 
$p_{m}^{\rm tree}(\beta)$ is given by
\begin{equation}\label{eq: inf}
p_{m}^{\rm tree}(\beta)=\inf_{\theta \in ]0,1]}v_{m}(\theta) 
\end{equation}
where $v_{m}$ is given by the expression
\begin{equation}\label{eq: exp}
\forall \theta \in ]0,1] \qquad
v_{m}(\theta)=\frac{1}{\theta}\ln(Q\sum_{x \in L_{m}}W_{m}(x)^{\theta}).
\end{equation}

We will first need the following monotonicity lemma.
\begin{lemma}\label{lem: mon}
Assume that $\phi: ]0,\infty[ \longrightarrow \ R$ is $\mathcal{C}^1$
and that there are constants $C,p \in [1, \infty[$ such that 
\begin{equation*}
\forall u > 0 \qquad \mid \phi'(u) \mid \leq C u^p + C u^{-p}.
\end{equation*}
Then for all $x \in L_{m}$ $\phi(W_{m}(x)),
\frac{\partial\phi(W_{m}(x))}{\partial\beta} \in L^{1}(Q)$,
$Q\phi(W_{m}(x))$ is $\mathcal{C}^1$ in $\beta \in \R$ and 
\begin{equation*}
\frac{\partial}{\partial\beta}Q\phi(W_{m}(x))=
Q\frac{\partial}{\partial\beta}\phi(W_{m}(x)).
\end{equation*}
Suppose in addition that $\phi$ is concave. Then ,
\begin{equation*} 
\forall \beta \geq 0 \qquad
Q\frac{\partial}{\partial\beta}\phi(W_{m}(x)) \leq 0.
\end{equation*}
\end{lemma}

\proof
The proof is an immediate adaptation of the proof of lemma 3.3 in
\cite{cf: CoYo}.
\qed

As a consequence we can define the following

\begin{proposition} \label{th:monotonicity}
The function $p_{m}^{\rm tree}$ is non-increasing in $\beta$.
There exists a critical value $\beta_{c}^{m} \in (0, \infty]$ such
that 
\begin{equation*}
p_{m}^{\rm tree}(\beta)=
\begin{cases}
0& \quad \text{if} \quad \beta \in [0,\beta_{c}^{m}], \\
< 0 & \quad \text{if} \quad \beta > \beta_{c}^{m}.
\end{cases}
\end{equation*}
\end{proposition}

\proof
For all $\theta \in ]0,1]$, the function $x \rightarrow x^{\theta}$ is
concave so by lemma \ref{lem: mon}, we see from expression (\ref{eq: exp}) that $v_{m}(\theta)$ is  non-increasing as a
function of $\beta$. Therefore, we see from (\ref{eq: inf}) that
$p_{m}^{\rm tree}$ is itself non-increasing in $\beta$ and we obtain
the existence of $\beta_{c}^{m}$ ($\beta_{c}^{m} \in [0, \infty]$). 
Since
$$
v_{m}'(1)
= Q\sum_{x \in L_{m}}W_{m}(x) \ln W_{m}(x)
\longrightarrow 
\sum_{x \in L_{m}} P(\omega_n=x) \ln  P(\omega_n=x) < 0\;,$$
as $\beta \searrow 0$, we conclude that  $\beta_{c}^{m}$ is strictly positive
by continuity of $\partial_\theta v_m (\theta, \beta)_{|\theta =1}$
in $\beta$ and by lemma \ref{lem: mon}. 

\qed

\begin{theorem}\label{lem: ine}
We have the following inequality
\begin{equation}\label{ine}
p(\beta) \leq \inf_{m \geq 1}\frac{1}{m} p_{m}^{\rm tree}(\beta).  
\end{equation}
\end{theorem}

\proof
Let $\theta \in (0,1)$ and $m$ be a positive integer. By using the
subadditive estimate 
\begin{equation}\label{eq: sub}
\forall u,v >0, \qquad (u+v)^{\theta} < u^{\theta}+v^{\theta}, 
\end{equation}
we have for all $n \geq 1$ 
\begin{align*}
Q \frac{1}{n}\ln W_{nm}& = 
Q \frac{1}{\theta n}\ln W_{nm}^{\theta} \\
& 
{\stackrel{(\ref{eq: mark})}{=}}
Q\frac{1}{\theta n}\ln\left(\sum_{x_{1}, \ldots,x_{n}}W_{m}(x_{1}) \ldots
 W_{(n-1)m,nm}^{x_{n-1}}(x_{n})\right)^{\theta} \\
& 
{\stackrel{(\ref{eq: sub})}{\leq}}
 Q\frac{1}{\theta n}\ln\sum_{x_{1}, \ldots,x_{n}}W_{m}(x_{1})^{\theta} \ldots
 W_{(n-1)m,nm}^{x_{n-1}}(x_{n})^{\theta} \\
& 
{\stackrel{\text{\tiny (Jensen)}}{\leq}}
 \frac{1}{\theta n}\ln Q \sum_{x_{1}, \ldots,x_{n}}W_{m}(x_{1})^{\theta} \ldots
 W_{(n-1)m,nm}^{x_{n-1}}(x_{n})^{\theta} \\
& =\frac{1}{\theta n}\ln\left(Q\sum_{x}W_{m}(x)^{\theta}\right)^{n} \\
& =\frac{1}{\theta}\ln Q\sum_{x}W_{m}(x)^{\theta} 
\end{align*}
The proof is complete by taking the limit as $ n \to \infty$ and then
by taking the infimum over all $\theta \in
]0,1]$ and $m \geq 1$.
\qed

In particular, to prove $p(\beta)<0$ it suffices to find $m \geq 1$ 
(in fact, $m \geq 2$) 
and $\theta \in (0,1)$ such that $Q \sum_{x}W_{m}(x)^{\theta}<1$. 
The theorem is a handy way to obtain upper bounds on the critical $\beta$. 

\begin{rem}\label{rem1}
Let $\theta \in ]0,1[$ and $m \geq 1$. Using (\ref{eq: sub}), we find by a similar computation that for
all $k \geq 2$ 
\begin{align}  \nonumber
Q\sum_{y}W_{km}(y)^{\theta} & =
Q\sum_{y}\big(\sum_{x_{1}, \ldots, x_{k-1}}W_{m}(x_{1}) \ldots W_{(k-1)m,km}^{x_{k-1}}(y)\big)^{\theta} \\ \nonumber
& < Q\sum_{y}\sum_{x_{1}, \ldots, x_{k-1}}W_{m}(x_{1})^{\theta} \ldots W_{(k-1)m,km}^{x_{k-1}}(y)^{\theta} \\ \label{eq: rem29}
& = \big(Q \sum_{x}W_{m}(x)^{\theta}\big)^k.
\end{align}
In view of (\ref{eq: inf}) and of the smoothness of $v_m(\cdot)$, 
we conclude that 
\begin{equation*}
\frac{1}{km} p_{km}^{\rm tree}(\beta) \leq \frac{1}{m} 
p_{m}^{\rm tree}(\beta). 
\end{equation*}
Observe that when $p_{m}^{\rm tree}(\beta)<0$, the infimum in  (\ref{eq: inf})
is achieved for some $\theta \in (0,1)$, and therefore the 
 above inequality is strict.
In particular,
\begin{equation} \label{eq: liminf-tree}
\inf_{m \geq 1}\frac{1}{m} p_{m}^{\rm tree}(\beta)=\varliminf_{m \to \infty}\frac{1}{m} p_{m}^{\rm tree}(\beta). 
\end{equation}
The authors do not know if the sequence 
$(p_{m}^{\rm tree}(\beta))_{m \geq 1}$ is  subadditive. However a simple 
argument yields the stronger result
\begin{equation}\label{eq: lim-tree}
\inf_{m \geq 1}\frac{1}{m} p_{m}^{\rm tree}(\beta)=\lim_{m \to \infty}\frac{1}{m} p_{m}^{\rm tree}(\beta). 
\end{equation}
Indeed, by repeating the steps in (\ref{eq: rem29}), we we see that, for 
$0 \leq \ell < m, k \geq 1$ and $\theta \in (0,1]$,
$$v_{km+\ell}(\theta) \leq k v_m(\theta) + v_\ell(\theta)\;,
$$
whereas, by concavity,
$$
 v_\ell(\theta) \leq \frac{1}{\theta} \sum_x \big(Q W_\ell(\theta)\big)^\theta
= v_\ell(\theta,0)
$$
where $v_\ell(\theta,0)= v_\ell(\theta, \beta)_{|\beta=0} \in (0, \8)$.
Therefore,
$$
\max_{km \leq n < (k+1)m} \frac{v_n(\theta)}{n} \leq
 \frac{k}{(k+\varepsilon)m} v_m(\theta) +  \frac{1}{km}v_\ell(\theta,0)\;,
$$
where $\varepsilon=0$ or $1$ according to the sign of $ v_m(\theta)$.
Now, recalling that $ v_m(\theta)\geq p_{m}^{\rm tree}(\beta)$ and taking 
the limit $k \to \8$, leads to
$$
\limsup_n \frac{p_{n}^{\rm tree}(\beta)}{n} \leq 
\frac{ v_m(\theta)}{m},\quad m \geq 1, \theta \in (0,1].
$$
Combined with (\ref{eq: liminf-tree}), this implies  (\ref{eq: lim-tree}).
\end{rem}
We add another
\begin{rem}
Suppose that there exists $m \geq 1$ such that
\begin{equation*} 
Q\sum_{x}W_{m}(x)\ln W_{m}(x)=0.
\end{equation*}
We have 
\begin{align*}
Q\sum_{y}W_{2m}(y)\ln W_{2m}(y) & =
Q\sum_{x,y}W_{m}(x)W_{m,2m}^{x}(y)
\ln W_{2m}(y) \\
& > \sum_{x,y}QW_{m}(x)W_{m,2m}^{x}(y)\ln
\left(W_{m}(x)W_{m,2m}^{x}(y)\right) \\
& = \sum_{x}\big(Q W_{m}(x)\ln W_{m}(x)\big)\sum_{y}QW_{m,2m}^{x}(y) \\
& \quad + \sum_{x}\big(QW_{m}(x)\big)
\sum_{y}Q W_{m,2m}^{x}(y)\ln W_{m,2m}^{x}(y) \\
& = 2\sum_{x}QW_{m}(x)\ln W_{m}(x) \\
& = 0
\end{align*}
Hence, by lemma \ref{lem: v}, $p_{2m}^{\rm tree}(\beta)<0$ and finally
$p(\beta)<0$.
\end{rem}
As a consequence of theorem \ref{lem: ine}, we get our main result
\medskip

\noindent
{\it Proof of theorem \ref{mainth}:} 
Let $\theta \in ]0,1]$ and $\beta > 0$. By using lemma 4.1 in \cite{cf: CoShYo}, there
exists a $c(\theta)>0$ such that 
\begin{equation*} 
\forall m \geq 1  \qquad  Q(W_{m}^{\theta}) \leq e^{-c(\theta)m^{\frac{1}{3}}}.\end{equation*}
Therefore 
\begin{align*}
Q(\sum_{x \in L_{m}}(W_{m}(x))^{\theta}) & \leq \mid L_{m} \mid
Q(W_{m}^{\theta}) \\
& \leq \mid L_{m} \mid e^{-c(\theta)m^{\frac{1}{3}}} \underset{m \to
  \infty}{\longrightarrow} 0,
\end{align*}
where we have used the fact that $\mid L_{m} \mid = O(m)$.
In particular, there exists $m \geq 1$ such that 
\begin{equation*}
Q(\sum_{x \in L_{m}}(W_{m}(x))^{\theta}) < 1.
\end{equation*}
We have $p_{m}^{\rm tree}(\beta) < 0$ and so by  theorem \ref{lem: ine}
 $p(\beta) < 0$.
\qed

\begin{theorem} \label{th:varfor}
Suppose the environment $\eta$ is bounded or gaussian. Then the
inequality (\ref{ine}) is in fact an equality
\begin{equation*} 
p(\beta)= \inf_{m \geq 1}p_{m}^{\rm tree}(\beta).
\end{equation*}
\end{theorem}
\proof
The inequality $p(\beta) \leq \inf_{m \geq 1}p_{m}^{\rm tree}(\beta)$
is in fact the conclusion of theorem \ref{lem: ine}
 and thus is true for all
environments.

We must show that $\inf_{m \geq 1}p_{m}^{\rm tree}(\beta) \leq p(\beta)$.
We treat the gaussian case, the bounded case being similar.
If $\beta \leq \beta_{c}$, we have by definition $p(\beta)=0$ and
since for all $m \geq 1$, $p_{m}^{\rm tree}(\beta) \leq 0$, the result
is obvious.
Suppose that $\beta$ is such that $\beta > \beta_{c}$. By definition of
$\beta_{c}$, $p(\beta) < 0$. Let $\theta \in ]0,1]$. We have by the
concentration result (\ref{eq: conc})
\begin{align*}
Q(W_{m}^{\theta}) & = e^{\theta Q(\ln(W_{m}))} Q(e^{\theta(\ln
  W_{m}-Q(\ln(W_{m}))} )\\ 
& \leq e^{\theta p(\beta)m+\frac{\beta^2\theta^2m}{2}}.
\end{align*}
For all $m \geq 1$, 
\begin{align*}
 \frac{1}{m} p_{m}^{\rm tree}(\beta) & \leq \frac{1}{\theta m}\ln(Q(\sum_{x \in L_{m}}(W_{m}(x))^{\theta}))\\ 
& \leq \frac{1}{\theta m} \ln(\mid L_{m} \mid) + \frac{1}{\theta
  m}\ln(Q(W_{m}^{\theta})) \\
& \leq \frac{1}{\theta m} \ln(\mid L_{m} \mid)+p(\beta) +
\frac{\beta^2\theta}{2} \\
& \underset{m \to \infty}{\longrightarrow}
p(\beta)+ \frac{\beta^2\theta}{2}
\end{align*}
where we have used the fact that $\mid L_{m} \mid =O(m^d)$. Thus, by
remark \ref{rem1}
\begin{equation*}
\inf_{m \geq 1}\frac{1}{m} p_{m}^{\rm tree}(\beta)=\lim_{m \to
  \infty}\frac{1}{m} p_{m}^{\rm tree}(\beta) \leq p(\beta)+ \frac{\beta^2\theta}{2}.
\end{equation*}
The proof is complete by letting $\theta \downarrow 0$. 
\qed

\end{document}